\documentclass[12pt, reqno]{amsart}
\usepackage{amssymb,amsthm,amsmath}
\usepackage[numbers,sort&compress]{natbib}
\usepackage{amssymb,amsmath}
\usepackage{amsfonts}
\usepackage{mathrsfs}
\usepackage{latexsym}
\usepackage{amssymb}
\usepackage{amsthm}
\usepackage{pdfsync}
\usepackage{indentfirst}
\usepackage{xcolor}
\date{Dec 26, 2019}

\newcommand{\R}{{\mathbb R}}

\renewcommand{\div}{{\rm div}}

\newcommand{\vphi}{\varphi}
\newcommand{\bd}{\bar\partial}

\let\oldsection\section
\renewcommand\section{\setcounter{equation}{0}\oldsection}

\newtheorem{corollary}{Corollary}[section]
\newtheorem{theorem}{Theorem}[section]

\newtheorem{remark}{Remark}[section]

\def\ba{\begin{eqnarray}}
\def\ea{\end{eqnarray}}

\def\R{\Bbb R}

\newcommand{\beq}{\begin{equation}}
\newcommand{\eeq}{\end{equation}}
\newcommand{\ben}{\begin{eqnarray}}
\newcommand{\een}{\end{eqnarray}}
\newcommand{\beno}{\begin{eqnarray*}}
\newcommand{\eeno}{\end{eqnarray*}}

\allowdisplaybreaks
\begin{document}

\title[Asymptotic behavior of the steady Navier-Stokes flow]{Asymptotic behavior of the steady Navier-Stokes flow in the exterior domain}

\author{Yueyang~Men}
\address[Yueyang~Men]{College of mathematics and physics, Beijing University of Chemical Technology\\ Beijing,
100029, China}
\email{menyuey@mail.buct.edu.cn}

\author{Wendong Wang}
\address[Wendong Wang]{School of Mathematical Sciences, Dalian University of Technology, Dalian, 116024,  China}
\email{wendong@dlut.edu.cn}

\author{Lingling Zhao}
\address[Lingling Zhao]{School of Mathematical Sciences, Dalian University of Technology, Dalian, 116024,  China}
\email{linglingz@126.com}



%%% ----------------------------------------------------------------------

\begin{abstract}
We consider an elliptic equation with unbounded drift in an exterior domain, and obtain quantitative uniqueness estimates at infinity, i.e. the non-trivial solution of $-\triangle u+W\cdot\nabla u=0$  decays in the form of $\exp(-C|x|\log^2|x|)$ at infinity provided $\|W\|_{L^\infty(\mathbb{R}^2\setminus B_1)}\lesssim 1$, which is sharp with the help of some counterexamples. These results also generalize the decay theorem by Kenig-Wang \cite{KW2015} in the whole space. As an application,
the asymptotic behavior of an incompressible fluid around a bounded obstacle is also considered. Specially for the two-dimensional case, we can improve the decay rate in \cite{KL2019} to $\exp(-C|x|\log^2|x|)$, where  the minimal decaying rate of $\exp(-C|x|^{\frac32+})$ is obtained  by Kow-Lin  in a recent paper \cite{KL2019} by using appropriate Carleman estimates.
\end{abstract}

%%% ----------------------------------------------------------------------
\maketitle

%\tableofcontents
\allowdisplaybreaks

\section{Introduction}\label{sec1}

In this note, we consider the steady Navier-Stokes equations in an exterior domain $\Omega$:
\begin{equation}\label{eq:NS-2D}
\left\{\begin{array}{llll}
-\mu\Delta v+v\cdot \nabla v+\nabla \pi=0,\quad {\rm in}\quad \Omega,\\
{\rm div }~ v=0,\quad {\rm in}\quad \Omega.
\end{array}\right.
\end{equation}
Without loss of generality, we assume that $B_1=\{x;|x|<1\}$ and $\Omega=\R^n\setminus{B_1}$.

A classical result of Finn (\cite{F65}), established in 1965, stated that if $n=3$, $v|_{\partial\Omega}=0$ and $v =o(1/|x|)$, then $v=0$. A few years later, in 1969, Dyer-Edmunds (\cite{DE1969}) showed that, if $v$ has bounded second derivatives and if
\begin{equation}\label{DE}
  v(x) = O(\exp(-\exp(\alpha|x|^3))),\ \ \  \text{for all}\ \ \  \alpha >0,
\end{equation}
then $v=0$. Note that Dyer-Edmunds' assumption (\ref{DE}) is stronger than Finn's assumption $v=o(1/|x|)$, but Dyer-Edmunds' result depends only on the local behavior of $v$ as $|x|\to \infty$.
In \cite{KL2019}, Kow-Lin show that the minimal decaying rate of any nontrivial solution $v$ is a bit greater than $\exp(-C|x|^{3/2+})$ at infinity in dimension $n\ge2.$
Note that the decay is far from the prior estimates. For example, Gilbarg-Weinberger \cite{GW1978} showed the velocity
 $v(x) = o(\log^{\frac12} |x|)$ and $|D v|\leq o(|x|^{-\frac34}(\log |x|)^{\frac98})$ provided that the Dirichlet energy is bounded in an exterior domain
For more references on this topic,  we refer to \cite{GNP1997,LUW2011,DI2017} and the references therein.

Next we focused on the two-dimensional case. Using the revised technique as in \cite{KSW2015} and suitable Caccioppoli-type inequality, we can improve the bound to $\exp(-C|x|\log^2|x|)$.
Our first result is as follows.
\begin{theorem}\label{thm:main}
Suppose that $v$ is a smooth solution of (\ref{eq:NS-2D}) with $\|v\|_{L^p(\Omega)}\lesssim 1$ with $2<p\leq \infty$. Then there exists a constant $C_0>0$ such that
\beno
M(R)=\inf_{|x_0|=R}\int_{B_1(x_0)}|v(y)|^2dy\geq \exp(-C_0 R\log^2(R)  )
\eeno
provided that
\beno
M(10)\neq 0.
\eeno

\end{theorem}

\begin{remark}
The above result improves the decay estimate in \cite{KL2019}, where they proved the expotential decay as $\exp(-C|x|^{3/2+})$.

%Note that the decay estimate here seems to be sharp, since for the entire plane, the similar estimate was proved in \cite{KW2015}.
\end{remark}

%We exploit the method by Kenig-Silvestre-Wang.In \cite{KSW2015}, the authors prove the quantitative version of Landis' conjecture in the plane with $V\ge 0$ .
For the exterior domain, we consider the second order elliptic equation with
drift term
\beq\label{driftequ}
-\triangle u+W\cdot\nabla u=0
\eeq
where $W = (W_1,W_2)$ is a real vector-valued function with $L^p$ bound for
$2<p\le\infty$. For 2D case, it's the vorticity of Navier-stokes equation if we denote $u={\rm curl} \ v$. Here we are interested in the lower bound of the decay rate for
any nontrivial solution $u$.

When $p = \infty$, the problem is related to Landis'
conjecture \cite{KL88}. That is, let $u$ be a solution of $-\Delta u+ Vu=0$ with $\|u\|_{L^{\infty}(\R^n)}\le C_0$ satisfying $|u(x)|\le C\exp(-C|x|^{1+})$, then $u\equiv 0$.
Landis' conjecture was disproved by Meshkov \cite{m92}, who constructed such $V$ and nontrivial $u$ satisfying $|u(x)|\le C\exp(-C|x|^{\frac 43})$. He also showed that if $|u(x)|\le C\exp(-C|x|^{\frac 43+})$, then $u\equiv 0$. (It should be noted that both $V$ and $u$ constructed by Meshkov are
complex-valued functions. It remains an open question whether Landis' conjecture is true
for real-valued $V$ and $u$ .)
But if we consider the equation (\ref{driftequ}) (or $-\Delta u+ W\cdot \nabla u+Vu=0$),
if one applies a suitable Carleman estimate to (\ref{driftequ}) and a scaling
devise in \cite{BK05}, the best exponent one can get is $2$, namely, under the same
conditions stated above except $|u(x)| \le \exp(-C|x|^{2+})$, then $u$ is trivial (see
\cite{Da12} for quantitative forms of this result). Moreover, in \cite{Da12}
the author constructed a Meshkov type example showing that the exponent
2 is in fact optimal for complex-valued $W$ and $u$ .
\par
In a recent paper \cite{KSW2015},  Kenig-Silvestre-Wang studied Landis' conjecture for
second order elliptic equations in the plane in the real setting, including (\ref{driftequ})
with real-valued $W$ and $u$. It was proved in \cite{KSW2015} that if $u$ is a real-valued
solution of $-\Delta u+ Vu=0$ with $V\ge 0$ satisfying $|u(x)|\le C\exp(-C|x|^{1+})$, then $u\equiv 0$. In fact, they prove the lower bound estimate for the nontrivial solution. For the equation (\ref{driftequ}) with $L^p$ drift term in entire plane, they prove the lower bound estimate depend on $p$. More references, we refer to
\cite{DKW2019,DZ2019} and the references therein.

In this paper, we prove the following decay rate in an exterior domain.
\begin{theorem}\label{thm:landis}
Let $u$ be a real solution of
\begin{equation}\label{eq:elliptic}
\Delta u-W\cdot \nabla u=0\quad\text{in}\quad B^c_1=\mathbb{R}^2\setminus{B_1},
\end{equation}
where $W$ satisfies
\begin{equation}\label{eq:W condi}
\|W\|_{L^{p}(B^c_1)}\le 1,\quad 2<p\leq \infty.
\end{equation}
Moreover, assume that $\|\nabla u\|_{L^{\infty}(B^c_1)}\leq 1$ and there exists $C_{0}>0$ such that
\begin{equation}\label{eq:lower bound}
\inf_{|z_{0}|=3}\int_{B_{1}(z_{0})}|\nabla u|^{2}\geq C_{0}.
\end{equation}
Then there holds
\ben\label{eq:decay}
\inf_{|z_{0}|=R}\sup_{|z-z_{0}|<1}|u(z)|\geq C_2\exp(-C_1R\log^2R),\ \  \text{for}\ \  ~~~~R\gg 1,
\een
where $C_2$ and $C_1$ are constants depending on $p$ and $C_0.$
\end{theorem}

%\wred{
\begin{remark}
Actually, the decay rates above are sharp when $p=\infty$. For example, $u(x)=\exp(-|x|)(1+|x|)$ and
$W(x)=\frac{2-|x|}{|x|^2}x\in L^{\infty}(B_1^c)$ solve the equation \eqref{eq:elliptic}.
\end{remark}

%\begin{remark}
%Moreover, by simple computations, one can verify that
%\beno
%u=\exp(-|x|^{1-\frac2p})(C_0+|x|)
%\eeno
%and
%\beno
%W(x)=\frac{(\frac4p-2)|x|^{-\frac2p}+|x|^{-1}+(C_0+|x|)(\frac2p-1)^2(|x|^{-\frac4p}-|x|^{-\frac2p-1})}{|x|+(C_0+|x|)(\frac2p-1)|x|^{-\frac2p+1}}x
%\eeno
%satisfy the equation $-\Delta u+W\cdot \nabla u=0$ in $B_1^c$ for a large $C_0>0$. Moreover, we also have
%\beno
%|W(x)|\sim |x|^{-\frac2p},\quad {\rm as}~|x|\rightarrow\infty,
%\eeno
%which implies that
%\beno
%W(x)\in L^{p+}(B_1^c).
%\eeno
%\end{remark}

%Consider the equation\\
% $-\Delta u+W\cdot \nabla u=0$ in $B_1^c$ in $\R^2$ with $W(x)=\frac{2-|x|}{|x|^2}\left(x_1,x_2\right)\in L^{\infty}(B_1^c)$ and
%\begin{equation*}
%   u(x)=e^{-|x|}(1+|x|)
%\end{equation*}
%is an exponent decay solution of the equation.

As in \cite{KW2015}, one can also replace the condition of \eqref{eq:lower bound} by the lower bound at every point.
\begin{corollary}\label{cor:landis2}
Let $u$ be a real solution of \eqref{eq:elliptic},
where $W$ satisfies
\eqref{eq:W condi}.
Moreover, assume that $\|\nabla u\|_{L^{\infty}(B^c_1)}\leq 1$ and there exists $C_{0}>0$ such that
\begin{equation}\label{eq:lower bound2}
\inf_{|z_{0}|=3}|\nabla u|\geq C_{0}'.
\end{equation}
Then there holds
\ben\label{eq:decay2}
\inf_{|z_{0}|=R}\sup_{|z-z_{0}|<1}|u(z)|\geq C_2'\exp(-C_1'R\log^2R),\ \  \text{for}\ \  ~~~~R\gg 1,
\een
where $C_2'$ and $C_1'$ are constants depending on $p$ and $C_0'.$
\end{corollary}

%We need the following lemmas.
Throughout the remaining section, we denote by $C$ a general positive constant which depends only on known constant coefficients and may be different from line to line.

\section{Proof of Theorem \ref{thm:main}}

In this section, we complete the proof of Theorem \ref{thm:main} with the help of Theorem \ref{thm:landis}.
One main step is to obtain the higher regularity of the equation (\ref{eq:NS-2D}) via the condition $\|v\|_{L^p(\Omega)}<\infty$.
We follow the same route as the proof of Liouville type theorems, for example see \cite{SW2019,WW2019}, where the divergence equation, Poincar\'{e}-Sobolev inequality and iteration lemma are used.

%
%First, we introduce a Giaquinta's iteration lemma \cite[Lemma 3.1]{G83}.
%
%\begin{lemma}[Lemma 3.1 \cite{G83}]\label{lem:Gia}
%\label{iter}
%%[\cite[Lemma 8]{ChYa}]
%Let $f(r)$ be a non-negative bounded function on $[R_0,R_1] \subset \R_+$. If there are negative constants $A,B,D$ and positive exponents $b<a$ and a parameter $\theta \in (0,1)$ such that for all $R_0 \le \rho < \tau \le R_1$
%\[
%f(\rho) \le \theta f(\tau) + \frac{A}{(\tau-\rho)^{a}} + \frac{B}{(\tau-\rho)^{b}} + D,
%\]
%then for all $R_0 \le \rho < \tau \le R_1$
%\[
%f(\rho) \leq C(a,\theta) \left[\frac{A}{(\tau-\rho)^{a}} + \frac{B}{(\tau-\rho)^{b}} + D\right].
%\]
%\end{lemma}

Next we begin to prove Theorem \ref{thm:main}.
\begin{proof}

{\bf Step I. Regularity estimates.} Assume that $B_{R}(x_0)\subset \Omega$ with $0<R\leq 1$, where $|x_0|\geq 2$, and
 $\mu=1$ without loss of generality. Choose a cut-off function $\phi(x)\in C_0^\infty(B_R(x_0))$ with $0\leq \phi\leq 1$ satisfying the following two properties:

\begin{enumerate}

\item[i).]
$\phi$ is radially decreasing and satisfies
\begin{align*} \phi(x)=\phi(|x-x_0|)=\left\{
\begin{aligned}
&1,\quad |x-x_0|\leq \rho,\\
&0, \quad |x-x_0|\geq\tau,
\end{aligned}
\right. \end{align*}
where $0<\frac{R}{2}\leq \rho<\tau\leq R$;

\item[ii).]
 $|\nabla\phi|(x)\leq \frac{C}{\tau-\rho}$, $|\nabla^2\phi|(x)\leq \frac{C}{(\tau-\rho)^2}$, $|\nabla^3\phi|(x)\leq \frac{C}{(\tau-\rho)^3}$ for all $x\in \mathbb{R}^2$.

\end{enumerate}

For $1<s<\infty$, due to Theorem III 3.1 in \cite{Galdi}, there exists a constant $C(s)$ and a vector-valued function $\bar{w}: B_\tau(x_0)\rightarrow \mathbb{R}^2$ such that
$\bar{w}\in W^{1,s}_0(B_\tau(x_0))$ and
$\nabla\cdot \bar{w}(x)=\nabla_x\cdot[\phi(x) {v}(x)]$. Moreover,  we get
\begin{align}
\label{esti-ws}
\int_{B_\tau(x_0)}|\nabla \bar{w}(x)|^s\,dx
\leq C(s)\int_{B_\tau(x_0)}|\nabla\phi\cdot { v}|^s\,dx.
\end{align}
%We thus can extend $\bar{w}$ to the whole space $\mathbb{R}^2$, which vanishes outside of the domain $B_\tau(x_0).$

Making the inner products $(\phi {v}-\bar{w})$ on both sides of the equation (\ref{eq:NS-2D}), by $\nabla\cdot \bar{w}=\nabla\cdot[\phi{v}]$ we have
\begin{align*}
&\hspace{-8mm}\int_{B_\tau(x_0)}\phi|\nabla v|^2\,dx
\\&= -\int_{B_\tau(x_0)}\nabla\phi\cdot\nabla v \cdot {v} \,dx+\int_{B_\tau(x_0)}\nabla \bar{w}:\nabla v  \,dx
-\int_{B_\tau(x_0)}v\cdot\nabla v \cdot \phi {v} \,dx\\
& \hspace{5mm}+\int_{B_\tau(x_0)}v\cdot\nabla v \cdot \bar{w} \,dx \\
&\doteq  I_1+\cdots+I_4,
\end{align*}
For the term $I_1$, it follows from H\"{o}lder inequality that
\beno
|I_1|\leq \frac{C}{\tau-\rho}\left(\int_{B_\tau(x_0)}|\nabla v|^2\,dx\right)^{\frac12}\left(\int_{B_\tau(x_0)}| {v}|^2\,dx\right)^{\frac12}.
\eeno
For the term $I_2$, H\"{o}lder inequality and (\ref{esti-ws}) imply that
\begin{align*}
|I_2|&\leq C\left(\int_{B_\tau(x_0)}|\nabla v|^2\,dx\right)^{\frac12}\|\nabla \bar{w}\|_{L^{2}(B_\tau(x_0))}\\
&\leq  \frac{C}{\tau-\rho}\|\nabla v\|_{L^2(B_\tau(x_0))} \|{v}\|_{L^2(B_\tau(x_0))}.
\end{align*}
By integration by parts and (\ref{esti-ws}), we find that
\beno
I_3+I_4\leq \frac{C}{\tau-\rho}\|v\|^{3}_{L^{3}(B_\tau(x_0))}.
\eeno
Combining the estimates of $I_1-I_4$,
\beno\label{eq:energy}
\int_{B_\tau(x_0)}\phi|\nabla v|^2\,dx\leq\frac14\|\nabla v\|_{L^2(B_\tau(x_0))}^2+ \frac{C}{(\tau-\rho)^{2}}\|v\|^{2}_{L^{2}(B_\tau(x_0))}+\frac{C}{\tau-\rho}\|v\|^{3}_{L^{3}(B_\tau(x_0))}.
\eeno
Recall that the following Poincar\'{e}-Sobolev inequality holds(see, for example, Theorem 8.11 and 8.12 \cite{LL})
\beno
\label{eq:poincare-sobolev}
\|f\|_{L^3(B_\tau)}\leq C \|\nabla f\|_{L^2(B_\tau)}^{\frac13}\|f\|_{L^2(B_\tau)}^{\frac23}+C\tau^{-\frac13}\|f\|_{L^2(B_\tau)},
\eeno
which implies that
\beno\label{eq:energy}
\int_{B_\tau(x_0)}\phi|\nabla v|^2\,dx&\leq&\frac12\|\nabla v\|_{L^2(B_\tau(x_0))}^2+ \frac{C}{(\tau-\rho)^{2}}\|v\|^{2}_{L^{2}(B_\tau(x_0))}\nonumber\\
&&+ \frac{C}{(\tau-\rho)^{2}}\|v\|^{4}_{L^{2}(B_\tau(x_0))}+\frac{C\tau^{-1}}{\tau-\rho}\|v\|^{3}_{L^{2}(B_\tau(x_0))}.
\eeno

Applying Giaquinta's iteration lemma (see \cite[Lemma 3.1]{G83}), we have
\ben\label{eq:L2 estimate}
\int_{B_\rho(x_0)}|\nabla v|^2\,dx\leq \frac{C}{(\tau-\rho)^{2}}\|v\|^{2}_{L^{2}(B_\tau(x_0))}+ \frac{C}{(\tau-\rho)^{2}}\|v\|^{4}_{L^{2}(B_\tau(x_0))}+\frac{C\tau^{-1}}{\tau-\rho}\|v\|^{3}_{L^{2}(B_\tau(x_0))}\nonumber\\
\een
Choose $\rho=R/2$ and $\tau=R$, and assume that $R=1$ without loss of generality. Since $\|v\|_{L^p(\Omega)}\lesssim 1$ with $2<p\leq \infty$,
\ben\label{eq:energy-1}
\int_{B_{1/2}(x_0)}|\nabla v|^2\,dx\leq C,
\een
for any $|x_0|\geq 2$.

Note that the vorticity $\omega=\partial_2v_1-\partial_1v_2$ is as follows:
\ben\label{eq:NS-2Dv}
-\Delta \omega+v\cdot \nabla \omega=0,\quad {\rm in}\quad \Omega.
\een
Making the inner products $\phi {\omega}$ on both sides of the equation (\ref{eq:NS-2Dv}), we have
\beno
\int_{B_R(x_0)}\phi|\nabla \omega|^2\,dx
&=& -\int_{B_R(x_0)}\nabla\phi\cdot\nabla \omega \cdot {\omega} \,dx
-\int_{B_R(x_0)}v\cdot\nabla \omega \cdot \phi {\omega} \,dx\\
&\doteq & I_1'+I_2'.
\eeno
For the term $I_1'$, by H\"{o}lder inequality  we have
\beno
|I_1'|\leq \frac{C}{\tau-\rho}\left(\int_{B_\tau(x_0)}|\nabla \omega|^2\,dx\right)^{\frac12}\left(\int_{B_R(x_0)}| {\omega}|^2\,dx\right)^{\frac12}\leq \frac{C}{\tau-\rho}\left(\int_{B_\tau(x_0)}|\nabla \omega|^2\,dx\right)^{\frac12},
\eeno
where we used (\ref{eq:energy-1}).
By integration by parts, we find that
\beno
I_2'=\int_{B_R(x_0)}v\cdot\nabla\phi \omega^2\,dx,
\eeno
Then
\beno
\int_{B_R(x_0)}\phi|\nabla \omega|^2\,dx
&\leq & \frac14\int_{B_\tau(x_0)}|\nabla \omega|^2\,dx+\frac{C}{(\tau-\rho)^2}+\frac{C}{\tau-\rho}\|v\|_{L^p(B_\tau(x_0))}\|\omega\|_{L^{2p'}(B_\tau(x_0))}^2
\eeno
where $\frac{1}{p'}+\frac1p=1.$
Note that when $p=\infty$, we have $p'=1$ and  $\|\omega\|_{L^{2p'}(B_\tau(x_0))}^2\leq C$ due to (\ref{eq:energy-1}).
Next, assume that $2<p<\infty.$
Using Poincar\'{e}-Sobolev inequality again
\beno
\label{eq:poincare-sobolev}
\|f\|_{L^{2p'}(B_\tau)}\leq C \|\nabla f\|_{L^2(B_\tau)}^{1-\frac1{p'}}\|f\|_{L^2(B_\tau)}^{\frac1{p'}}+C\tau^{-1+\frac1{p'}}\|f\|_{L^2(B_\tau)},
\eeno
which implies that
\beno
\int_{B_R(x_0)}\phi|\nabla \omega|^2\,dx
&\leq & \frac12\int_{B_\tau(x_0)}|\nabla \omega|^2\,dx+\frac{C}{(\tau-\rho)^2}+\frac{C}{(\tau-\rho)^{p'}}+\frac{C}{\tau-\rho}\tau^{-2+\frac2{p'}}
\eeno
Applying Giaquinta's iteration lemma again, we have
\ben\label{eq:energy-2}
\int_{B_{1/2}(x_0)}|\nabla \omega|^2\,dx\leq C,
\een
for any $|x_0|\geq 3$. In fact, \eqref{eq:energy-2} implies that
\ben\label{eq:energy-3}
\int_{B_{1/2}(x_0)}|\nabla^2 v|^2\,dx\leq C,
\een
due to integration by parts and $\triangle v=\nabla\div(v)-{\rm curl}{\rm curl} v.$ Moreover, (\ref{eq:energy-3}) and $\|v\|_{L^p(\Omega)}\lesssim 1$
yields that
\ben\label{eq:v bounded}
\|v\|_{L^\infty(\mathbb{R}^2\setminus{B_3})}\leq C,
\een
by Gagliardo-Nirenberg inequality.

Furthermore, using the equation \eqref{eq:NS-2Dv} we get
\beno
\int_{B_{1}(x_0)}|\triangle \omega|^2\,dx\leq C\int_{B_{1}(x_0)}|v|^2|\nabla \omega|^2\,dx\leq C,
\eeno
where we used \eqref{eq:v bounded} and \eqref{eq:energy-3}. It follows that
\ben\label{eq:energy-nabla3}
\int_{B_{1}(x_0)}|\nabla^3v|^2\,dx\leq C,
\een
which and (\ref{eq:energy-1}) yield that
\ben\label{eq:nabla v bound}
\|\nabla v\|_{L^\infty(\mathbb{R}^2\setminus{B_3})}\leq C.
\een

Similarly, using the equation \eqref{eq:NS-2Dv} again,
\beno
\int_{B_{1}(x_0)}|\triangle\nabla \omega|^2\,dx\leq C\int_{B_{1}(x_0)}|\nabla (v\cdot\nabla \omega)|^2\,dx\leq C,
\eeno
where we used \eqref{eq:v bounded}, \eqref{eq:energy-nabla3}, (\ref{eq:energy-2}), \eqref{eq:nabla v bound},  and Gagliardo-Nirenberg inequality. It follows that
\beno\label{eq:nabla4}
\int_{B_{1}(x_0)}|\nabla^4 v|^2\,dx\leq C,
\eeno
which and \eqref{eq:energy-3} yield that
\ben\label{eq:nabla v bound'}
\|\nabla^2 v\|_{L^\infty(\mathbb{R}^2\setminus{B_3})}\leq C.
\een

{\bf Step II. Decay estimates of the vorticity.}

Note that the vorticity satisfies the maximum principle, then there exist constants $C_{0}'$ and $R_0>2$ such that
\beno
\inf_{|x_1|=R_0}\int_{B_1(x_1)}|\nabla \omega|^2dx\geq C_{0}',
\eeno
since $M(10)\neq 0.$

Applying Theorem \ref{thm:landis} due to \eqref{eq:nabla v bound'} and scaling property, by (\ref{eq:decay}) we have
\ben\label{eq:decay2'}
\inf_{|x_{0}|=R}\int_{B_1(x_0)}|\omega|^2dx\geq C_2'\exp(-C_1'R\log^2(R)),\ \  \text{for}\ \  ~~~~R\gg 1.
\een

{\bf Step III. Decay estimates of the velocity.}

By the energy inequality \eqref{eq:L2 estimate} and (\ref{eq:v bounded}), we have
\beno
\inf_{|x_{0}|=R}\int_{B_1(x_0)}|\omega|^2&\leq& \inf_{|x_{0}|=R}\int_{B_1(x_0)}| \nabla v|^2dx\\&
\leq&\inf_{|x_{0}|=R}\int_{B_1(x_0)}| v|^2dx.
\eeno
which and \eqref{eq:decay2'} imply that
\beno
M(R)=\inf_{|x_0|=R}\int_{B_1(x_0)}|v(y)|^2dy\geq \exp(-C_0 R\log^2(R)  )
\eeno
Then the proof is complete.

\end{proof}

\section{Proof of Theorem \ref{thm:landis}}

\begin{proof} We follow the same route as in \cite{KSW2015}. The difference is,  we choose a different cut-off function due to the exterior domain and deal with the $L^p$ drift.

Let $z_0'\in\R^2$ with $|z_0'|\gg 1$. Since \eqref{eq:elliptic} is invariant under rotation, we can assume that $z'_0=|z_0'|e_1$, where $e_1=(1,0)$. Translating the origin to $-3e_1$, \eqref{eq:elliptic} becomes
\begin{equation}\label{eq:u exterior}
\Delta u-W(x,y)\cdot\nabla u=0\quad\text{in}\quad B^c_1(-3e_1).
\end{equation}
For simplicity, we still write $u$ and $W$ in the equation in the new coordinates. Now we denote $z_0=(|z_0'|-3)e_1$ and set $R=|z_0|$.  Define the scaled solution $u_R(z)=u(ARz+z_0)$, where $A>0$, to be decided. Therefore, $u_R$ solves
\begin{equation}\label{eq:u-R}
\Delta u_R-W_R\cdot\nabla u_R=0\quad\text{in}\quad B_{\frac{1}{AR}}^c(z_1),
\end{equation}
where
\[
z_1=-(\frac 1A+\frac{3}{AR})e_1
\]
and $W_R(z)=(AR)W(ARz+z_0)$. Thus, for any $2<p\leq \infty$ there holds
\ben\label{eq:w-R}
\|W_R\|_{L^{p}(B_{\frac{1}{AR}}^c(z_1))}\le (AR)^{1-\frac2p},
\een
where we used (\ref{eq:W condi}).
And the origin ($ARz+z_0=0$) moves to
\[
\hat z=-\frac{z_0}{AR}=-\frac 1Ae_1.
\]
Choose a large $A$ so that
\[
B_{\frac{1}{AR}}(z_1)\subset B_{7/5}.
\]

Note that $\triangle=4\partial\bar{\partial}$, where
\beno
\partial=\frac12(\partial_x-i\partial_y),\quad \bar{\partial}=\frac12(\partial_x+i\partial_y).
\eeno
It follows from (\ref{eq:u-R}) that $u_R$ satisfies
\beno
4\partial\bar{\partial}u_R-W_R\cdot\left((\partial+\bar{\partial})u_R, -i(\bar{\partial}-\partial)u_R\right)=0,
\eeno
which implies
\beno
\bar{\partial}(\partial u_R)=\alpha\partial u_R,
\eeno
where we define
\ben\label{eq:alpha}
\alpha\doteq\frac14 W_R\cdot\left(1+\frac{\bar{\partial} u_R}{\partial u_R}, -i\frac{\bar{\partial} u_R}{\partial u_R}+i\right),
\een
for $|z-z_1|\geq \frac{1}{AR}$, otherwise $\alpha=0.$

Let
$g=\chi\partial u_R$, here $\chi$ is a cutoff function $\chi\equiv 1$ on $|z-z_1|\ge \frac{9}{8AR}$ and $\chi\equiv 0$ for $|z-z_1|\le\frac{17}{16AR}$. Note that $\nabla\chi$  is supported on $\frac{17}{16AR}\le|z-z_1|\le\frac{9}{8AR}$.
Then we have
\begin{equation}\label{eq: g' equ}
\bar\partial g=\alpha g+\bar{\partial}\chi\partial u_R \quad\text{in}\quad B_2.
\end{equation}
We now write $\hat z$ as a point in the complex plane, i.e., $\hat z=-\frac{1}{A}+i0$. Let $w(z)$ be defined by
\[
w(z)=\frac{1}{\pi}\int_{B_{7/5}}\frac{\alpha}{\xi-z}d\xi-\frac{1}{\pi}\int_{B_{7/5}}\frac{\alpha}{\xi-\hat z}d\xi,
\]
then $\bar{\partial} w=-\alpha$ in $B_{7/5}$. Recalling that (\ref{eq:w-R}) and (\ref{eq:alpha}), we have
\beno
\|\alpha\|_{L^{p}(B_{7/5})}\le C(AR)^{1-\frac2p},
\eeno
In view of \cite{ve62} (for example, see (6.4)-(6.7),(6.9a)), we have the following estimate of $w(z)$.
For $2<p<\infty$, there holds
\ben\label{eq:2<p<infty}
|w(z)|\le C(p) \|\alpha\|_{L^{p}(B_{7/5})}|z-\hat z|^{1-\frac2p}\leq C(AR)^{1-\frac2p}|z-\hat z|^{1-\frac2p},\quad\forall\ z\in B_{7/5},
\een
and for $p=\infty$
\ben\label{eq:infty}
|w(z)|\le C(AR)|z-\hat z|\log\left(\frac{C}{|z-\hat z|}\right),\quad\forall\ z\in B_{7/5}.
\een

Let $h=e^wg$, then it follows from \eqref{eq: g' equ} that
\ben\label{eq: h'}
\bar{\partial} h=e^{w}(\bd\chi)\partial u_R \quad\text{in}\quad B_{7/5},
\een
Next we will use the following Carleman type estimate of $\bd$ form \cite[Proposition~2.1]{df90}. Let $\vphi_{\tau}(z)=\vphi_{\tau}(|z|)=-\tau\log|z|+|z|^2$, then for any $f\in C_0^{\infty}(B_{7/5}\setminus\{0\})$, we have that
\begin{equation}\label{care}
\int |\bd f|^2e^{\vphi_{\tau}}\ge\frac 14\int(\Delta\vphi_{\tau})|f|^2e^{\vphi_{\tau}}=\int|f|^2e^{\vphi_{\tau}}.
\end{equation}
Note that $\vphi_{\tau}$ is decreasing in $|z|$ for $\tau>8$ and $|z|\leq 2$. We introduce another cutoff function $0\le\zeta\le 1$ satisfying
\[
\zeta(z)=\left\{
\begin{aligned}
0,&\quad\text{when}\ |z|<\frac{1}{4AR},\\
1,&\quad\text{when}\ \frac{1}{2AR}<|z|<1,\\
0,&\quad\text{when}\ |z|>6/5.
\end{aligned}\right.
\]
Hence the following estimates holds
\beno\label{estz}
|\nabla\zeta(z)|\le C(AR)\;\;\text{for}\;\; z\in X\;\;\text{and}\;\;|\nabla\zeta(z)|\le C\;\;\text{for}\;\; z\in Y,
\eeno
where
\[
X=\{\frac{1}{4AR}<|z|<\frac{1}{2AR}\}\;\; \text{and}\;\; Y=\{1<|z|<6/5\}.
\]
We also denote
\[
Z=\{\frac{1}{2AR}<|z|<1\}.
\]
Note that (\ref{eq: h'}), and applying the Carleman estimate \eqref{care} to $\zeta h$ we have
\begin{equation}\label{eq:h' estimate}
\begin{aligned}
\int_Z|h|^2e^{\vphi_{\tau}}&\le 2\int(|\bd\zeta h|^2+|\zeta\bd h|^2)e^{\vphi_{\tau}}\\
&\le C(AR)^2\int_X|h|^2e^{\vphi_{\tau}}+C\int_Y|h|^2e^{\vphi_{\tau}}+\int_{\widetilde Z}|e^{w}(\bd\chi)\partial u_R|^2e^{\vphi_{\tau}},
\end{aligned}
\end{equation}
where
\[
\widetilde Z=\{\frac{1}{4AR}<|z|<\frac 65\}.
\]
First, for the left one of  (\ref{eq:h' estimate}), for $A$ and $R$ large enough we have
\[
\int_Z|h|^2e^{\vphi_{\tau}}\ge\int_{B_{\frac{1}{AR}}(\hat z)}|h|^2e^{\vphi_{\tau}}.
\]
Next let us estimate the terms in the integral inequality of (\ref{eq:h' estimate}) in two cases.

{\bf Step I. Case I of  $2<p<\infty$.} On one hand,  it follows from \eqref{eq:2<p<infty} that
\beno
|w(z)|\le C;\;\text{for}\;\; z\in B_{\frac{1}{AR}}(\hat z),
\eeno
i.e.,
\beno\label{eq:w lower bound p>2}
e^{w(z)}\ge\frac{1}{C}\;\;\text{for}\;\; z\in B_{\frac{1}{AR}}(\hat z).
\eeno
And using that for $z\in B_{\frac{1}{AR}}(\hat z)$, $|z|\le\frac{1}{AR}+\frac 1A$,  we have
\begin{equation}\label{bar}
\int_{Z}|h|^2e^{\vphi_{\tau}}\ge\frac{e^{\vphi_{\tau}(\frac 1A+\frac{1}{AR})}}{C}\int_{B_{\frac{1}{AR}}(\hat z)}|\partial u_R|^2.
\end{equation}
Next we look at $\int_{\widetilde Z}|e^{w}(\bd\chi)\partial u_R|^2e^{\vphi_{\tau}}$. Recall that $\bd\chi$ is supported in
 $\frac{17}{16AR}\le|z-z_1|\le\frac{9}{8AR}$. Thus
\begin{equation}\label{ewz1}
e^{w(z)}\le C\;\;\text{for}\;\; \frac{17}{16AR}\le|z-z_1|\le\frac{9}{8AR}.
\end{equation}
Using \eqref{ewz1} and the known condition $\|\nabla u\|_\infty\lesssim 1$, we have
\begin{equation}\label{h1e}
\int_{\widetilde Z}|e^{w}(\bd\chi)\partial u_R|^2e^{\vphi_{\tau}}\le C(AR)^2\int_{\frac{17}{16AR}\le|z-z_1|\le\frac{9}{8AR}}|\partial u_R|^2e^{\vphi_{\tau}}\le C(AR)^2e^{\vphi_{\tau}(\frac 1A+\frac{15}{8AR})}.
\end{equation}

It follows from \eqref{eq:2<p<infty} that
\beno
|w(z)|\leq C(AR)^{1-\frac2p},\quad\forall\ z\in B_{7/5}.
\eeno
Multiplying $\exp(-\vphi_{\tau}(\frac 1A+\frac{1}{AR}))$ on both sides of \eqref{eq:h' estimate}, using \eqref{bar}, \eqref{h1e} and the  bound of $\partial u$, we obtain
\begin{equation}\label{h12e}
\begin{aligned}
\int_{B_{\frac{1}{AR}(\hat z)}}|\partial u_R|^2&\le C(AR)^2e^{C(AR)^{1-\frac2p}}\frac{\exp(\vphi_{\tau}(\frac{1}{4AR}))}{\exp(\vphi_{\tau}(\frac 1A+\frac{1}{AR}))}\int_{B_{\frac{1}{2AR}(0)}}|\partial u_R|^2\\
&\;\;\;+Ce^{C(AR)^{1-\frac2p}}\frac{\exp(\vphi_{\tau}(1))}{\exp(\vphi_{\tau}(\frac 1A+\frac{1}{AR}))}\\
&\;\;\;+C(AR)^2\frac{\exp({\vphi_{\tau}(\frac 1A+\frac{15}{8AR})})}{\exp(\vphi_{\tau}(\frac 1A+\frac{1}{AR}))}\\
\end{aligned}
\end{equation}

Re-scaling back to the original variables, by (\ref{eq:lower bound}) we observe that
\ben\label{eq:rescaling}
\int_{B_{\frac{1}{AR}(\hat z)}}|\partial u_R|^2=\int_{B_1(0)}|\partial u|^2\ge {C_0}\;\;\quad\text{and}\quad\int_{B_{\frac{1}{2AR}(0)}}|\partial u_R|^2=\int_{B_{\frac12}(z_0)}|\partial u|^2
\een
Finally, choosing $\tau=C(AR)\log(AR)$ and taking $R$ sufficiently large, it is not hard to see that
\[
\left\{\begin{aligned}
&C(AR)^{2}e^{C(AR)^{1-\frac2p}}\frac{\exp(\vphi_{\tau}(\frac{1}{4AR}))}{\exp(\vphi_{\tau}(\frac 1A+\frac{1}{AR}))}\le\exp(C(AR)\log^2(AR)),\\
&C(AR)^{2}e^{C(AR)^{1-\frac2p}}\frac{\exp(\vphi_{\tau}(1))}{\exp(\vphi_{\tau}(\frac 1A+\frac{1}{AR}))}\to 0,\\
&C(AR)^{2}\frac{\exp({\vphi_{\tau}(\frac 1A+\frac{15}{8AR})})}{\exp(\vphi_{\tau}(\frac 1A+\frac{1}{AR}))}\to 0,\\
\end{aligned}\right.
\]
Therefore, if $R$ is large enough, then the last two term on the right hand side of \eqref{h12e} can be absorbed by the term on the left. Hence we get
\begin{equation}\label{patialv}
  \int_{B_{\frac12}(z_0)}|\partial u|^2\ge C\exp(-C(AR)\log^2(AR)).
\end{equation}
Note that the energy estimate implies
\beno
\int_{B_{\frac12}(z_0)}|\partial u|^2\leq C \int_{B_1(z_0)}| u|^2+C\left(\int_{B_1(z_0)}| u|^2\right)^{1/2}
\eeno
which yields that
\beno
\int_{B_1(z_0)}| u|^2\geq C \min\left\{\int_{B_{\frac12}(z_0)}|\partial u|^2, \left(\int_{B_{\frac12}(z_0)}|\partial u|^2\right)^{2}\right\}
\eeno
Then the case of $2<p<\infty$ is complete.

At last, we deal with the case of $p=\infty.$

{\bf Step II. Case II. $p=\infty$.} On the other hand,
\eqref{eq:infty} implies
\[
|w(z)|\le C\ln(AR);\;\text{for}\;\; z\in B_{\frac{1}{AR}}(\hat z),
\]
hence
\beno\label{eq:w lower bound infty}
e^{w(z)}\ge\frac{1}{(AR)^C}\;\;\text{for}\;\; z\in B_{\frac{1}{AR}}(\hat z).
\eeno
Similarly as \eqref{bar},  we have
\begin{equation}\label{bar'}
\int_{Z}|h|^2e^{\vphi_{\tau}}\ge\frac{e^{\vphi_{\tau}(\frac 1A+\frac{1}{AR})}}{(AR)^C}\int_{B_{\frac{1}{AR}}(\hat z)}|\partial u_R|^2.
\end{equation}
and
\begin{equation}\label{h1e'}
\int_{\widetilde Z}|e^{w}(\bd\chi)\partial u|^2e^{\vphi_{\tau}}\le C(AR)^C\int_{\frac{17}{16AR}\le|z-z_1|\le\frac{9}{8AR}}|\partial u_R|^2e^{\vphi_{\tau}}\le C(AR)^Ce^{\vphi_{\tau}(\frac 1A+\frac{15}{8AR})}.
\end{equation}
Note that
\beno
|w(z)|\leq C(AR),\quad\forall\ z\in B_{7/5}.
\eeno
Multiplying $C(AR)^C\exp(-\vphi_{\tau}(\frac 1A+\frac{1}{AR}))$ on both sides of \eqref{eq:h' estimate}, using \eqref{bar'}, \eqref{h1e'} and the  bound of $\partial u$, we obtain
\begin{equation}\label{h12e'}
\begin{aligned}
\int_{B_{\frac{1}{AR}(\hat z)}}|\partial u_R|^2&\le C(AR)^C\exp(C(AR))\frac{\exp(\vphi_{\tau}(\frac{1}{4AR}))}{\exp(\vphi_{\tau}(\frac 1A+\frac{1}{AR}))}\int_{B_{\frac{1}{2AR}(0)}}|\partial u_R|^2\\
&\;\;\;+C(AR)^C\exp(C(AR))\frac{\exp(\vphi_{\tau}(1))}{\exp(\vphi_{\tau}(\frac 1A+\frac{1}{AR}))}\\
&\;\;\;+C(AR)^C\frac{\exp({\vphi_{\tau}(\frac 1A+\frac{15}{8AR})})}{\exp(\vphi_{\tau}(\frac 1A+\frac{1}{AR}))}\\
\end{aligned}
\end{equation}

Re-scaling back to the original variables again as in \eqref{eq:rescaling}, there hold
\[
\int_{B_{\frac{1}{AR}(\hat z)}}|\partial u_R|^2=\int_{B_1(0)}|\partial u|^2\ge C_0,\quad\text{and}\quad\int_{B_{\frac{1}{2AR}(0)}}|\partial u_R|^2=\int_{B_{\frac12}(z_0)}|\partial u|^2
\]
Finally, choosing $\tau=C(AR)\log(AR)$ and taking $R$ sufficiently large, we have
\[
\left\{\begin{aligned}
&C(AR)^{C}\exp(C(AR))\frac{\exp(\vphi_{\tau}(\frac{1}{4AR}))}{\exp(\vphi_{\tau}(\frac 1A+\frac{1}{AR}))}\le\exp(CAR(\log^2(AR))),\\
&C(AR)^{C}\exp(C(AR))\frac{\exp(\vphi_{\tau}(1))}{\exp(\vphi_{\tau}(\frac 1A+\frac{1}{AR}))}\to 0,\\
&C(AR)^{C}\frac{\exp({\vphi_{\tau}(\frac 1A+\frac{15}{8AR})})}{\exp(\vphi_{\tau}(\frac 1A+\frac{1}{AR}))}\to 0,\\
\end{aligned}\right.
\]
Therefore, if $R$ is large enough, then the last two term on the right hand side of \eqref{h12e'} can be absorbed by the term on the left. Consequently,
\begin{equation}\label{patialv'}
  \int_{B_{\frac12}(z_0)}|\partial u|^2\ge C \exp(-CAR(\log^2(AR))).
\end{equation}
Hence the proof is complete by the interior estimate as in Step I.

\end{proof}

\section{Proof of Corollary \ref{cor:landis2}}

\begin{proof}

Since $W\in L^p$ with $p>2$ and $\|\nabla u\|_\infty\leq 1$, we have
\beno
u\in W^{2,p}_{loc}(\mathbb{R}^2\setminus B_1),
\eeno
which implies that
$u\in C^1(\mathbb{R}^2\setminus B_1)$.
Furthermore, by \eqref{eq:lower bound2}
there exists a positive constant $\delta$ such that
\beno
\inf_{|z_{0}|=3}\int_{B_{\delta}(z_{0})}|\nabla u|^{2}\geq C_{0}.
\eeno
With the help of Theorem \ref{thm:landis}, the proof is complete.

\end{proof}

\section*{Acknowledgments}
{ W. Wang was supported by NSFC under grant 11671067 and
 "the Fundamental Research Funds for the Central Universities".
}
\par

\end{document}